\documentclass[12pt,reqno]{amsart}

\newcommand\version{September 19, 2011}

\usepackage{amsmath,amsfonts,amsthm,amssymb,amsxtra}

\newtheorem{theorem}{Theorem}[section]

\newtheorem{lemma}[theorem]{Lemma}

\theoremstyle{definition}

\theoremstyle{remark}

\numberwithin{equation}{section}

\newcommand{\C}{\mathbb{C}}

\renewcommand{\epsilon}{\varepsilon}

\newcommand{\N}{\mathbb{N}}

\renewcommand{\phi}{\varphi}
\newcommand{\R}{\mathbb{R}}

\newcommand{\Sph}{\mathbb{S}}

\DeclareMathOperator{\spa}{span}

\begin{document}

\title[Ground states of non-local equations --- \version]{On the uniqueness of
ground states\\ of non-local equations}

\author{Rupert L. Frank}
\address{Rupert L. Frank, Department of Mathematics,
Princeton University, Washington Road, Princeton, NJ 08544, USA}
\email{rlfrank@math.princeton.edu}

\thanks{\copyright\, 2011 by the author. This paper may be reproduced, in its
entirety, for non-commercial purposes.\\
This work was supported through the U.S.~National Science Foundation grant
PHY-1068285.}

\begin{abstract}
We review our joint result with E. Lenzmann about the uniqueness of ground state
solutions of non-linear equations involving the fractional Laplacian and
provide an alternate uniqueness proof for an equation related to the
intermediate long-wave equation.
\end{abstract}

\maketitle

\section{Introduction}

We are interested in uniqueness of positive solutions $Q$ of non-linear
equations, of which
\begin{equation}\label{eq:eqintro}
(-\Delta)^s Q - Q^{\alpha+1} = -Q
\quad\text{in}\ \R^N
\end{equation}
with $0<s<1$ is a typical example. For $s$ in this range, $(-\Delta)^s$ is a
non-local operator. Throughout we will only be interested in finite energy
solutions, i.e., $Q\in H^s(\R^N)$.

\medskip

Before discussing \eqref{eq:eqintro} in more detail, we recall a few classical
results about the local case with $s=1$,
\begin{equation}\label{eq:eqintro1}
-\Delta Q - Q^{\alpha+1} = -Q
\quad\text{in}\ \R^N \,.
\end{equation}
These results concern
\emph{existence}, \emph{radial symmetry} and \emph{uniqueness} of positive
solutions. The existence of non-negative solutions for sub-critical values
of $\alpha$ (i.e., $2<\alpha+2<\infty$ if $N=1,2$ and $2<\alpha+2<2^*=2N/(N-2)$
if $N\geq 3$) follows from the fact that the infimum
$$
\inf_{u\in H^1(\R^N)} \frac{\left(\int |\nabla u|^2 \,dx \right)^\theta
\left(\int |u|^2 \,dx \right)^{1-\theta}}{\left(\int |u|^{\alpha+2} \,dx
\right)^{2/(\alpha+2)} } \,,
\qquad \theta= \frac{N\alpha}{2(\alpha+2)} \,,
$$
is attained. Any minimizer, after scaling and multiplication by a constant,
yields a non-negative solution to \eqref{eq:eqintro1}. By rearrangement
inequalities this shows that \eqref{eq:eqintro1} has a positive solution
which is symmetric decreasing (i.e., radial with respect to some point and
non-increasing with respect to the distance from this point). It is a theorem
that \emph{any} positive solution is (strictly) symmetric decreasing. This can
be proved by the method of moving planes. The most delicate result about
\eqref{eq:eqintro1} is due to Kwong \cite{Kw} (extending earlier results in
\cite{Co,MLSe}) 
and states uniqueness of positive solutions vanishing at infinity. The proof is
based on a careful ODE analysis. To summarize, for any subcritical $\alpha$,
\eqref{eq:eqintro1} has a unique positive, finite-energy solution.

\medskip

The question we would like to address here is whether the same is true for the
non-local equation
\eqref{eq:eqintro} with $0<s<1$. Equations of this type arise in numerous
models from mathematical physics, mathematical biology and mathematical finance
and we refer to \cite{FrLe2} for precise references. Although
\eqref{eq:eqintro} is time-independent, it plays a role in the description of
traveling wave or solitary wave solutions of certain time-dependent equations.

Existence of positive solutions follows again by minimization for sub-critical
values of $\alpha$. Now `sub-critical' means $2<\alpha+2<\infty$ if $N\leq 2s$
and $2<\alpha+2<2N/(N-2s)$ if $N>2s$. Radial symmetry and monotonicity of
positive solutions has recently been established in \cite{MaZh} (see also
\cite{FrLe1}) by a modification of the method of moving planes. Thus, in what
follows we shall concentrate on the \emph{uniqueness} of positive, radial
solutions of \eqref{eq:eqintro}.

To understand why a uniqueness result is plausible and where the assumption
$0<s<1$ comes from, we recall that the heat kernel of $(-\Delta)^s$, that is,
the
integral kernel of the operator $\exp(-t(-\Delta)^s)$, is positive precisely
for $0<s\leq 1$. By well-known arguments, this implies that positive solutions
$Q$ of the linear Schr\"odinger-type equation $(-\Delta)^s Q + V Q = -\mu Q$ are
unique (up to a multiplicative constant). This suggests that the
uniqueness property is closer related to the positivity of the heat kernel (and
therefore to the maximum principle) than it is to locality and ODE techniques,
which are crucial in Kwong's work \cite{Kw}.

The question of uniqueness has been raised for ground states of some specific
non-local equations (see, e.g., \cite{LYa} and the preprint version of
\cite{KeMaRo}), but rigorous answers are rare in the literature. A celebrated
result of Amick and Toland \cite{AmTo} concerns the case $N=1$, $s=1/2$ and
$\alpha=1$ of \eqref{eq:eqintro}. Later, we will discuss in more detail a
related result of Albert and Toland \cite{AlTo}, again for $N=1$ and $\alpha=1$,
but with a non-homogeneous operator. The only other case we are aware of
concerns the critical case $\alpha+2=2N/(N-2s)$ of \eqref{eq:eqintro} with the
right side replaced by zero. Here \emph{any} $0<s<N/2$ is allowed, see
\cite{L,ChLiOu}. We emphasize that in all these cases the optimizers are known
explicitly and the uniqueness proofs make use of these closed form expressions.

\medskip

Recently, in \cite{FrLe2} we proved a uniqueness result for energy-minimizing
positive solutions of \eqref{eq:eqintro} in dimension $N=1$. The parameter
$\alpha$ is allowed to take any sub-critical value. The precise statement is the
following.

\begin{theorem}\label{fl}
 Let $N=1$ and assume that $0<\alpha<\frac{4s}{1-2s}$ if $0<s<\frac12$ and that
$0<\alpha<\infty$ if $\frac12\leq s<1$. Then the infimum
$$
\inf_{u\in H^s(\R^N)} \ \frac{\|(-\Delta)^{s/2} u \|_2^{\theta} \
\|u\|_2^{1-\theta}}{\|u\|_{\alpha+2}}\,,
\qquad \theta= \frac{\alpha}{2s(\alpha+2)} \,,
$$
is attained by a unique (modulo symmetries) function. More precisely,
there is a positive, even and decreasing function $Q$ such that any
minimizer is of the form $a Q(b(x-c))$ for some $a\in \C\setminus\{0\}$, $b>0$
and $c\in\R$.

Moreover, $Q$ satisfies $(-\Delta)^s Q - Q^{\alpha+1}=-Q$, and the corresponding
linearization
$$
L_+ = (-\Delta)^s -(\alpha+1)Q^{\alpha}+1
\qquad\mathrm{in}\ L^2(\R)
$$
is non-degenerate, i.e., $\ker L_+ = \spa \{Q'\}$.
\end{theorem}

\noindent\emph{Remarks.}
(1) The restriction on $\alpha$ are optimal for the existence of a minimizer.\\
(2) Most (but not all) of the arguments in our proof work for $N\geq 2$ as
well.\\
(3) It is an \emph{open problem} whether every positive solution is energy
minimizing.\\
(4) No closed form expression for the functions $Q$ are known unless $s=1/2$ and
$\alpha=1$. Their regularity and (inverse power) decay are rather well
understood.\\
(5) The fact that $Q'\in\ker L_+$ follows by differentiating the equation for
$Q$. The difficult part is to prove that there is no linear independent element
in $\ker L_+$. This so-called non-degeneracy result has applications in
stability and blow-up analysis of non-linear dispersive equations; see, e.g.,
\cite{KeMaRo}.

\medskip

The idea behind our proof of Theorem \ref{fl} is to find a continuous branch of functions $Q_s$, $s\in[s_0,1]$, such that $Q_s$ solves \eqref{eq:eqintro} for any $s$ and such that $Q_{s_0}$ agrees with a given energy-minimizing solution of the equation. Uniqueness at $s_0$ is then deduced from uniqueness and non-degeneracy at $s=1$, which are well known. The construction of the branch is based on the inverse function theorem, which requires one to check that a certain linear operator is invertible. To prove this invertibility we need a slightly stronger non-degeneracy statement than that in Theorem \ref{fl}. A crucial ingredient of
its proof is a result about the number of zeroes of an eigenfunction corresponding to the second eigenvalue of a Schr\"odinger-type operator $(-\Delta)^s +V$. It is only at this last point that we need the assumption that
$N=1$. For the details of the proof we refer the reader to \cite{FrLe2}.

\subsection*{Acknowledgment}
I wish to thank E. Lenzmann in collaboration with whom Theorem \ref{fl} was
obtained. His questions stimulated me to find the results presented in the
following sections.


\section{A uniqueness result for the ILW equation}

We turn now to the uniqueness question for a one-dimensional equation which is
slightly different from \eqref{eq:eqintro}. It arises in the analysis of the
intermediate long-wave (ILW) equation, which interpolates, in some sense,
between the well-known KdV and Benjamin--Ono equations; see \cite{Jo} and the
references therein for background information. Travelling wave solutions of the
ILW equation correspond to solutions $Q\in H^{1/2}(\R)$ of the equation
\begin{equation}
 \label{eq:ilw}
T Q - Q^2 = -\mu Q
\quad\text{in}\ \R \,.
\end{equation}
Here $T$ is the (non-local) pseudo-differential operator given by
$$
\widehat{ T\phi}(\tau) = \left( \tau \coth(\pi\tau/2)-2/\pi \right)
\widehat\phi(\tau)
$$
in terms of the Fourier transform $\widehat\phi(\tau)=(2\pi)^{-1/2} \int
e^{-i\tau t}\phi(t)\,dt$, and we are using dimensionless
variables in which the height $H=\pi/2$. Our results concern the special case
$\mu=2/\pi$ of \eqref{eq:ilw}. We shall prove

\begin{theorem}\label{ineq}
 For any $\phi\in H^{1/2}(\R)$ one has
$$
(\phi,T\phi) + \tfrac2\pi \|\phi\|_2^2 \geq \left(\tfrac \pi2 \right)^{1/3}
\|\phi\|_3^2
$$
with equality iff $\phi(x)= \beta\, h(t-\alpha)$ for some $\alpha\in\R$ and
$\beta\in\C$, where
$$
h(t) =\frac1{\cosh t} \,.
$$
\end{theorem}

This result implies the existence of a positive solution to equation
\eqref{eq:ilw} with $\mu=2/\pi$ and states that it is unique (up to
translations) among energy-minimizing solutions. The next result shows that,
actually, the energy-minimizing property is not necessary for this uniqueness
result. In particular, any positive solution is automatically energy-minimizing.

\begin{theorem}\label{eq}
 Assume $Q\in H^{1/2}(\R)$ is non-negative and satisfies
\begin{equation}
 \label{eq:eq}
T Q - Q^2 = - \tfrac2\pi Q 
\qquad\text{in}\ \R \,.
\end{equation}
Then $Q(x)= h(t-\alpha)$ for some $\alpha\in\R$. Moreover,
$$
\ker\left( T-2Q+ \tfrac2\pi \right) = \spa\{Q'\} \,.
$$
\end{theorem}

\noindent\emph{Remarks.}
(1) The first part of this theorem is due to Albert and Toland \cite{AlTo} (see
also \cite{Al}) and the second part due to Albert and Bona \cite{AlBo}. They
can treat equation \eqref{eq:ilw} with general $\mu>0$.\\
(2) The assumption $Q\geq 0$ in Theorem \ref{eq} can be replaced by the weaker
assumption that $Q$ is real-valued. To see this, write the equation as
$Q=(T+\tfrac2\pi)^{-1} Q^2$ and use the fact that $(T+\mu)^{-1}$ has a positive
integral
kernel for any $\mu>0$; see, e.g., \cite{Al}.\\
(3) The fact that $Q'\in\ker\left( T-2Q+\tfrac2\pi \right)$ follows from
translation-invariance of \eqref{eq:eq}. The non-obvious fact is that there is
no linearly independent function in this kernel.\\
(4) Equation \eqref{eq:eq} is understood in $H^{-1/2}(\R)$, that is, for any
$\phi\in H^{1/2}(\R)$
$$
\int_\R \tau \coth(\pi\tau/2)\, \overline{\widehat\phi(\tau)} \widehat Q(\tau)
\,d\tau = \int_\R \overline{\phi(t)} Q(t)^2 \,dt \,.
$$

\medskip

The reason why we think our new proofs of Theorems \ref{ineq} and \ref{eq} may
be worth recording is that they establish a connection between equation
\eqref{eq:eq} and the seemingly unrelated equation
$$
\sqrt{-\Delta} R= R^2
\quad\text{in}\ \R^3 \,.
$$
This will allow us to deduce Theorems \ref{ineq} and \ref{eq} from known results
about the latter equation. Here is our key lemma.

\begin{lemma}\label{key}
 Let $\phi$ be a function on $\R$ and $\psi$ be a radial function on $\R^3$
related by
$$
\psi(x) = |x|^{-1} \phi(\ln |x|) \,.
$$
Then $\phi\in H^{1/2}(\R)$ iff $\psi\in\dot H^{1/2}(\R^3)$ and, in this case,
\begin{equation*}
 (\psi,\sqrt{-\Delta}\psi) = 4\pi \int_\R \tau \coth(\pi\tau/2)
|\widehat\phi(\tau)|^2 \,d\tau \,.
\end{equation*}
\end{lemma}

By way of comparison, we note that if $\psi(x) = |x|^{-1/2} \phi(\ln |x|)$, then
$$
\int_{\R^3} |\nabla \psi|^2 \,dx = 4\pi \int_\R \left( |\phi'|^2 + \tfrac14
|\phi|^2\right) \,dt \,.
$$

Actually, Lemma \ref{key} works both ways: Not only can we deduce Theorems
\ref{ineq} and \ref{eq} from three-dimensional results, we also obtain sharp
functional inequalities in three dimensions from the one-dimensional
Albert--Toland result \cite{AlTo}.

\begin{theorem}\label{katosob}
 For any $\psi\in H^{1/2}(\R^3)$ and any $\theta\in(0,\pi/2)$ one has
$$
(\psi,\sqrt{-\Delta}\psi) - \tfrac2\pi \theta\cot\theta\ (\psi,|x|^{-1}\psi)  
\geq \left( \frac{4\theta}{\sqrt\pi \sin\theta} \right)^{\frac23} \left(
\theta\left(2+\cos(2\theta)\right) - \tfrac32 \sin(2\theta) \right)^{\frac13}\
\|\psi\|_3^2 \,,
$$
with equality iff $\psi(x) = c H_\theta(b x)$ for some $b>0$ and $c\in\C$, where
\begin{equation*}
H_\theta(x) = |x|^{-1} \left( |x|^{2\theta/\pi} + |x|^{-2\theta/\pi} +
2\cos\theta \right)^{-1} \,.
\end{equation*}
\end{theorem}

We note that $\theta\cot\theta$ strictly decreases as function of $\theta$ from
$1$ at $\theta=0$ to $0$ at $\theta=\pi/2$. Thus Theorem \ref{katosob}
interpolates between Kato's inequality $(\psi,\sqrt{-\Delta}\psi) \geq (2/\pi)
(\psi,|x|^{-1}\psi)$ and the Sobolev inequality $(\psi,\sqrt{-\Delta}\psi) \geq
2^{1/3} \pi^{2/3} \|\psi\|_3^2$ (see Lemma \ref{ineq} below).


\section{Proofs}

\subsection{Proof of the key lemma}

The ground state substitution formula from \cite{FrLiSe} (see also \cite{FrSe})
reads
$$
(\psi,\sqrt{-\Delta}\psi) - \frac{2}{\pi} \int_{\R^3} \frac{|\psi(x)|^2}{|x|} dx
= \frac{1}{2\pi^2} \iint_{\R^3\times\R^3}
\frac{|\chi(x)-\chi(y)|^2}{|x-y|^4}\frac{dx}{|x|} \frac{dy}{|y|} \,,
$$
where $\psi(x)= |x|^{-1} \chi(x)$. This is true for any $\psi\in\dot
H^{1/2}(\R^3)$. Now assume that $\psi$ is radial. Then, using the fact that for
$|x|=r$
$$
\int_{\Sph^2} \frac{d\omega}{|x-s\omega|^4} = 2\pi \int_0^\pi \frac{\sin\theta
d\theta}{(r^2-2rs\cos\theta+s^2)^2}
= 2\pi \int_{-1}^1 \frac{dt}{(r^2-2rst+s^2)^2} = \frac{4\pi}{(r^2-s^2)^2}\,,
$$
we obtain
$$
\iint_{\R^3\times\R^3} \frac{|\chi(x)-\chi(y)|^2}{|x-y|^4}\frac{dx}{|x|}
\frac{dy}{|y|}
= (4\pi)^2 \iint_{\R_+\times\R_+} \frac{|\chi(r)-\chi(s)|^2}{(r^2-s^2)^2} r
s\,dr\,ds \,.
$$
Next, we change variables $r=e^t$ and $s=e^u$. Writing $\chi(r)= \phi(\ln r)$ we
find
\begin{align*}
\iint_{\R_+\times\R_+} \frac{|\chi(r)-\chi(s)|^2}{(r^2-s^2)^2} r s\,dr\,ds
& = \iint_{\R\times\R} \left| \phi(t) - \phi(u)\right|^2
\frac{e^{2(t+u)}}{(e^{2t}-e^{2u})^2} \,dt\,du \\
& = \frac14 \iint_{\R\times\R} \frac{| \phi(t) - \phi(u)|^2}{\sinh^2(t-u)}
\,dt\,du \,.
\end{align*}
Summarizing what we have shown so far is that
$$
(\psi,\sqrt{-\Delta}\psi) = 2 \iint_{\R\times\R} \frac{| \phi(t) -
\phi(u)|^2}{\sinh^2(t-u)} \,dt\,du
+ 8 \int_{\R} |\phi(t)|^2 \,dt \,.
$$
Now we take the Fourier transform and find
\begin{align}
 \label{eq:form}
\iint_{\R\times\R} \frac{| \phi(t) - \phi(u)|^2}{\sinh^2(t-u)} \,dt\,du
& = 4 \int_\R \left( \int_\R \frac{\sin^2(\tau t/2)}{\sinh^2 t} \,dt \right)
|\widehat\phi(\tau)|^2 \,d\tau \notag \\
& = 8\pi \int_\R \left( \frac{\tau \coth(\pi\tau/2)}{4} - \frac{1}{2\pi} \right)
|\widehat\phi(\tau)|^2 \,d\tau \,,
\end{align}
where the last equality used \cite[(3.986.4)]{GrRy}. This proves Lemma
\ref{key}. \qed


\subsection{Proofs of Theorems \ref{ineq} and \ref{eq}}

Given our Key Lemma \ref{key}, we can deduce Theorems \ref{ineq} and \ref{eq}
from two known results in $\R^3$. The first is the following sharp Sobolev
inequality due, in a dual form, to Lieb \cite{Li}; see also \cite[Thm.
8.4]{LiLo}.

\begin{theorem}\label{lieb}
 For any $\psi\in \dot H^{1/2}(\R^3)$ one has
\begin{equation}
 \label{eq:lieb}
(\psi,\sqrt{-\Delta}\psi) \geq 2^{1/3} \pi^{2/3} \|\psi\|_3^2 \,,
\end{equation}
with equality iff $\psi(x) = c H(b(x-a))$ for some $a\in\R^3$, $b>0$ and
$c\in\C$, where
\begin{equation}
 \label{eq:h}
H(x) = (1+ |x|^2)^{-1} \,.
\end{equation}
\end{theorem}

If we combine this theorem with Lemma \ref{key}, we obtain
\begin{align*}
& 4\pi \int_\R \tau \coth(\pi\tau/2) |\widehat\phi(\tau)|^2 \,d\tau
= (\psi,\sqrt{-\Delta}\psi) \\ 
& \qquad \geq 2^{1/3} \pi^{2/3} \left( \int_{\R^3} |\psi|^3 \,dx \right)^{2/3}
= 2^{5/3} \pi^{4/3} \left( \int_\R |\phi|^3 \,dt \right)^{2/3} \,,
\end{align*}
with equality iff $\psi(x) = c\, H(b x)$ for some $b>0$ and $c\in\C$, that is,
$$
\phi(t) = \frac{c}{2b\, \cosh(t+\ln b)}
$$
This proves Theorem \ref{ineq}.

\medskip

For the proof of Theorem \ref{eq} we need the following result of Li \cite{L}
and Chen--Li--Ou \cite{ChLiOu}, which had been conjectured by Lieb \cite{Li}.
Since we have not been able to find a proof of the non-degeneracy statement in
the literature we provide one in the appendix.

\begin{theorem}\label{lclo}
 Assume that $R\in \dot H^{1/2}(\R^3)$ is non-negative and satisfies
\begin{equation}
 \label{eq:lclo}
\sqrt{-\Delta} R = R^2
\qquad\text{in}\ \R^3 \,.
\end{equation}
Then  $R(x) = 2b\, H(b(x-a))$ for some $a\in\R^3$ and $b\geq 0$ with $H$ from
\eqref{eq:h}. Moreover,
\begin{equation}
 \label{eq:nondeg3}
\ker\left( \sqrt{-\Delta} -2R \right) = \spa\{\partial_1 R,\partial_2
R,\partial_3 R, R+x\cdot\nabla R\} \,.
\end{equation}
\end{theorem}

The proof of Theorem \ref{eq} is now similar to that of Theorem \ref{ineq}.
Indeed, polarization of the identity in Lemma \ref{key} shows that equation
\eqref{eq:eq} for $Q$ becomes equation \eqref{eq:lclo} for
$R(x)=|x|^{-1}Q(\ln|x|)$. Thus, by Theorem \ref{lclo}, $R(x)=2b H(bx)$, that is,
$Q(t) = \left(\cosh(t+\ln b)\right)^{-1}$, as claimed.

To prove the non-degeneracy statement for $Q$ we use the fact that, since $R$ is
centered at the origin, the operator $\sqrt{-\Delta} -2R$ commutes with
rotations about the origin. Therefore the subspace of radial functions in
$L^2(\R^3)$ is a reducing subspace for this operator and, according to
\eqref{eq:nondeg3}, its kernel restricted to this subspace is spanned by
$R+x\cdot\nabla R$. On the other hand, Lemma \ref{key} implies that there is a
one-to-one correspondence between the kernel of $T-2Q+\frac2\pi$ and the kernel
of $\sqrt{-\Delta} -2R$ on radial functions. The assertion now follows from the
fact that $R+x\cdot\nabla R=|x|^{-1}Q'(\ln|x|)$ if $R(x)=|x|^{-1}Q(\ln|x|)$.
This completes the proof of Theorem \ref{eq}.


\subsection{Proof of Theorem \ref{katosob}}

By rearrangement inequalities \cite[Thm. 3.4 and Lemma 7.17]{LiLo} (which are
\emph{strict} for the $(\psi,|x|^{-1}\psi)$-term), we only need to consider
radial functions. On such functions, the theorem is equivalent to
\begin{equation}
 \label{eq:ilwineq}
(\phi,T\phi) + \tfrac2\pi\left(1-\theta\cot\theta\right) \|\phi\|_2^2 
\geq \left( \frac{2\theta}{\pi \sin\theta} \right)^{2/3} \left(
\theta\left(2+\cos(2\theta)\right) - \tfrac32 \sin(2\theta) \right)^{1/3}
\|\phi\|_3^2 \,,
\end{equation}
with equality iff $\phi(t)= \beta h_\theta(t-\alpha)$ for some $\alpha\in\R$ and
$\beta\in\C$, where
$$
h_\theta(t) = \frac{(2\theta/\pi)\sin\theta}{\cosh(2\theta t/\pi) + \cos\theta}
\,.
$$
This follows from \cite{AlTo}. Indeed, standard arguments (see, e.g.,
\cite{AlBoSa}) yield that the infimum
$$
I_\theta = \inf_{\phi\in H^{1/2}(\R)} \frac{(\phi,T\phi) +
\tfrac2\pi\left(1-\theta\cot\theta\right) \|\phi\|_2^2}{\|\phi\|_3^2}
$$
is attained by a function $\phi_\theta$. Using \eqref{eq:form} one can show that
any minimizer is a multiple of a non-negative function, and therefore we may
assume that $\phi_\theta\geq 0$. The Euler--Lagrange equation reads
$$
T\phi_\theta + \tfrac2\pi\left(1-\theta\cot\theta\right) \phi_\theta = \lambda
\phi_\theta^2
$$
with a Lagrange multiplier $\lambda$. Integration against $\phi_\theta$ shows
that $I_\theta = \lambda \|\phi_\theta\|_3$. The function $Q=\lambda
\phi_\theta$ satisfies \eqref{eq:ilw} and therefore, by \cite{AlTo},
$Q(t)=h_\theta(t)$. Finally, by a computation,
\begin{align*}
 I_\theta & = \|Q\|_3 = \frac{2\theta\sin\theta}{\pi} \left( 2 \int_0^\infty
\frac{dt}{\left(\cosh(2\theta t/\pi) + \cos\theta \right)^3} \right)^{1/3} \\
& = \left( \frac{2\theta}{\pi \sin\theta} \right)^{2/3} \left(
\theta\left(2+\cos(2\theta)\right) - \tfrac32 \sin(2\theta) \right)^{1/3} \,,
\end{align*}
which concludes the proof of Theorem \ref{katosob}.

\medskip

\noindent\emph{Remark.} Inequality \eqref{eq:ilwineq} is true for all
$\theta\in(0,\pi)$ (by the same proof) and therefore the inequality in Theorem
\ref{katosob} is true in this range \emph{for radial functions} $\psi$. However,
rearrangement inequalities for the $(\psi,|x|^{-1}\psi)$-term work in the wrong
way for $\theta\in(\pi/2,\pi)$.


\appendix

\section{Non-degeneracy for $\sqrt{-\Delta}$ in three dimensions}

In this appendix we prove \eqref{eq:nondeg3}. Our argument is based on Lieb's
observation \cite{L} that inequality \eqref{eq:lieb}, and therefore also
equation \eqref{eq:lclo}, are conformally invariant. To utilize this fact, we
first turn the differential equation into an integral equation by noting that
$$
\ker\left( \sqrt{-\Delta} -2R \right) = \ker\left( 2 R^{1/2} (-\Delta)^{-1/2}
R^{1/2} - 1 \right) \,.
$$
Indeed, $\sqrt{-\Delta}\psi =2R\psi$ iff $2 R^{1/2} (-\Delta)^{-1/2} R^{1/2}\chi
=\chi$, where $\chi=R^{1/2}\psi$. (In mathematical physics, this is called the
\emph{Birman--Schwinger principle}.) The integral kernel of the operator
$R^{1/2} (-\Delta)^{-1/2} R^{1/2}$ appearing on the right side is
\begin{equation}
 \label{eq:kernel}
(2\pi^2)^{-1} \sqrt{R(x)} |x-y|^{-2} \sqrt{R(y)} \,.
\end{equation}

To proceed, we note that we may assume that $a=0$ and $b=1$ in the formula for
$R$. We introduce the stereographic projection $\mathcal S: \R^3\to\Sph^3 =
\{(\omega_1,\ldots,\omega_4)\in\R^4:\ \sum \omega_j^2 = 1\}$,
\begin{align*}
 \mathcal S_j(x) = \frac{2x_j}{1+x^2}\quad\text{for}\ j=1,\ldots,3\,,
 \qquad
\mathcal S_{4}(x) = \frac{1-x^2}{1+x^2} \,.
\end{align*}
The Jacobian of $\mathcal S$ is
$$
\mathcal J(x) = \left( \frac{2}{1+x^2} \right)^3 \,,
$$
and therefore, we obtain a unitary operator $\mathcal U:L^2(\Sph^3) \to
L^2(\R^3)$ by setting
\begin{equation}
\label{eq:fu}
\mathcal U \Psi(x) = \mathcal J(x)^{1/2} \Psi(\mathcal S(x)) \,,
\end{equation}
for $\Psi\in L^2(\Sph^3)$. Using the fact that
$$
|\mathcal S(x) - \mathcal S(y)|^2 = \mathcal J(x) |x-y|^2 \mathcal J(y) \,,
$$
we find that the operator $K=\mathcal U^* R^{1/2} (-\Delta)^{-1/2} R^{1/2}
\mathcal U$ has kernel
$$
K(\omega,\eta) = (2\pi^2)^{-1} |\omega-\eta|^{-2} \,.
$$
Since this operator commutes with rotations of $\Sph^3$, its eigenfunctions are
spherical harmonics and its eigenvalues can be computed explicitly using the
Funk--Hecke formula, see \cite{Be} and also \cite{FrLi}. The eigenvalues are
given by $(l+1)^{-1}$, $l\in\N_0$, and the multiplicity of the $l$-th
eigenvalue is $(l+1)^2$. Here $l$ is the degree of the spherical harmonic. We
conclude that
$$
\ker\mathcal U^* \left( 2 R^{1/2} (-\Delta)^{-1/2} R^{1/2} - 1 \right) \mathcal
U
= \ker\left( 2K-1 \right)
$$
corresponds to $l=1$ and is therefore spanned by the functions
$\Psi_j(\omega)=\omega_j$, $j=1,\ldots,4$. Thus $\ker\left( 2 R^{1/2}
(-\Delta)^{-1/2} R^{1/2} -
1 \right)$ is spanned by the $\,\mathcal U \Psi_j = \mathcal J^{1/2} \mathcal
S_j
=R^{3/2}\mathcal S_j$ and, again by the Birman--Schwinger principle,
$\ker\left( \sqrt{-\Delta} -2R \right)$ is spanned by the $R \mathcal S_j$.
These are the claimed functions and our proof of \eqref{eq:nondeg3} is complete.


\bibliographystyle{amsalpha}

\end{document}